\documentclass[11pt,A4paper]{article}

\usepackage[left=30mm,right=30mm,top=25mm,bottom=30mm]{geometry}
\usepackage{labelfig}
\usepackage{epsfig}
\usepackage{epstopdf}
\usepackage{soul}
\usepackage{xfrac}

\usepackage[percent]{overpic}

\usepackage{color}
\usepackage{amsthm,amsmath,amssymb}
\usepackage{booktabs}
\usepackage{mathpazo}
\usepackage{microtype}
\usepackage{overpic}
\usepackage{bm}
\usepackage{sectsty}
\usepackage[
	pdftitle={PDFTitle},
	pdfauthor={Hugo Parlier},
	ocgcolorlinks,
	linkcolor=linkred,
	citecolor=linkred,
	urlcolor=linkblue]
{hyperref}

\definecolor{linkred}{RGB}{0,191,255} 
\definecolor{linkblue}{RGB}{16, 78, 139}

\usepackage[hang,flushmargin]{footmisc}
\usepackage{enumitem}
\usepackage{titlesec}
	\titlespacing{\section}{0pt}{12pt}{0pt}
	\titlespacing{\subsection}{0pt}{6pt}{0pt}
	
\titlelabel{\thetitle.\quad}

\makeatletter 

\long\def\@footnotetext#1{%
\H@@footnotetext{%
\ifHy@nesting 
\hyper@@anchor{\@currentHref}{#1}%
\else 
\Hy@raisedlink{\hyper@@anchor{\@currentHref}{\relax}}#1%
\fi 
}}

\def\@footnotemark{%
\leavevmode 
\ifhmode\edef\@x@sf{\the\spacefactor}\nobreak\fi 
\H@refstepcounter{Hfootnote}%
\hyper@makecurrent{Hfootnote}%
\hyper@linkstart{link}{\@currentHref}%
\@makefnmark 
\hyper@linkend 
\ifhmode\spacefactor\@x@sf\fi 
\relax 
}%

\ifFN@multiplefootnote%
\renewcommand*\@footnotemark{%
\leavevmode 
\ifhmode 
\edef\@x@sf{\the\spacefactor}%
\FN@mf@check 
\nobreak 
\fi 
\H@refstepcounter{Hfootnote}%
\hyper@makecurrent{Hfootnote}%
\hyper@linkstart{link}{\@currentHref}%
\@makefnmark 
\hyper@linkend 
\ifFN@pp@towrite 
\FN@pp@writetemp 
\FN@pp@towritefalse 
\fi 
\FN@mf@prepare 
\ifhmode\spacefactor\@x@sf\fi 
\relax%
}%
\fi 

\makeatother 

\theoremstyle{plain}
\newtheorem{theorem}{Theorem}[section]
\newtheorem{proposition}[theorem]{Proposition}
\newtheorem{lemma}[theorem]{Lemma}
\newtheorem{corollary}[theorem]{Corollary}

\theoremstyle{definition}

\newcommand{\sys}{{\rm sys}}

\newcommand{\N}{{\mathbb N}}

\newcommand{\arcsinh}{{\,\rm arcsinh}}
\newcommand{\arccosh}{{\,\rm arccosh}}

\newcommand{\injrad}{\mathrm{injrad}}

\sectionfont{\large \bfseries}
\subsectionfont{\normalsize}

\long\def\symbolfootnote[#1]#2{\begingroup%
\def\thefootnote{\fnsymbol{footnote}}\footnote[#1]{#2}\endgroup}

\def\blfootnote{\xdef\@thefnmark{}\@footnotetext}

\usepackage{mathtools}
\DeclarePairedDelimiter\floor{\lfloor}{\rfloor}
\DeclarePairedDelimiter\ceil{\lceil}{\rceil}
\usepackage[sort,nocompress]{cite}

\begin{document}

{\Large \bfseries The shortest non-simple closed geodesics on hyperbolic surfaces}

{\large Ara Basmajian\symbolfootnote[1]{\small 
Supported by a grant from the Simons foundation (359956, A.B.) and PSC-CUNY Award 65245-00 53.},
Hugo Parlier\symbolfootnote[2]{\small Supported by the Luxembourg National Research Fund OPEN grant O19/13865598.
}
and Hanh Vo}

\vspace{0.5cm}
{\bf Abstract.}
This article explores closed geodesics on hyperbolic surfaces. We show that, for sufficiently large $k$, the shortest closed geodesics with at least $k$ self-intersections, taken among all hyperbolic surfaces, all lie on an ideal pair of pants and have length $2\arccosh(2k+1)$. 
\vspace{0.5cm}

\section{Introduction} \label{sec:intro}
The study of closed geodesics on hyperbolic surfaces has multiple facets which links together topics as diverse as spectral theory, symbolic dynamics, geometric topology and, of course, hyperbolic geometry. Without going into its rich history, in the past few decades, simple closed geodesics have played a capital role in the understanding of Teichm\"uller and moduli spaces, whereas closed geodesics were more traditionally related to spectral theory. The study of closed geodesics, and in particular the non-simple ones, has seen a recent resurgence over the last decade and a plethora of recent results suggest that stratifying closed geodesics by their self-intersection number, and asking similar questions to the ones that were asked for simple curves, provides fascinating directions of research. 

One way to do this is to study the family of curves with fixed self-intersection. Notable examples of successes in that direction include generalizations of Mirzakhani's simple curve counting theorem \cite{Mirzakhani2008} to curves with (fixed) self-intersection (see for instance \cite{Mirzakhani2016, Erlandsson-Souto, Erlandsson-Parlier-Souto, Sapir}). 

In a somewhat orthogonal direction, there has been a lot of work relating how intersection and length are related, in particular when the intersection grows. For each surface and each natural number, there is a shortest closed geodesic with (at least) that number of self-intersections. These are a sort of generalization of the systole (the shortest closed geodesic on a surface), and their geometry and actual number of self-intersections was first systematically studied in \cite{Basmajian1993}. While a lot of work has gone into extending these results \cite{Basmajian2013, Basmajian2015, Chas2015, Erlandsson-Parlier}, to this point it has been unclear whether these shortest curves have exactly the minimum number of self-intersections. The most precise result to date \cite{Vo} says that if a surface has a cusp, then for large enough $k$, the shortest geodesic with at least $k$ self-intersections has exactly $k$ self-intersections. As we will see, this result will provide a crucial step in our investigation. 

Let $\gamma_k$ be the closed geodesic on a connected hyperbolic surface $X_k$ of minimal length among all primitive closed geodesics with at least $k$ self-intersection points taken among all closed geodesics on all possible hyperbolic surfaces. Our main result is the following. 

\begin{theorem}\label{thm:main}
For large enough $k$, $X_k$ is an ideal pair of pants and $\gamma_k$ is a corkscrew geodesic with exactly $k$ self-intersections.
\end{theorem}

\begin{figure}[h]
\begin{center}
\includegraphics[width=6cm]{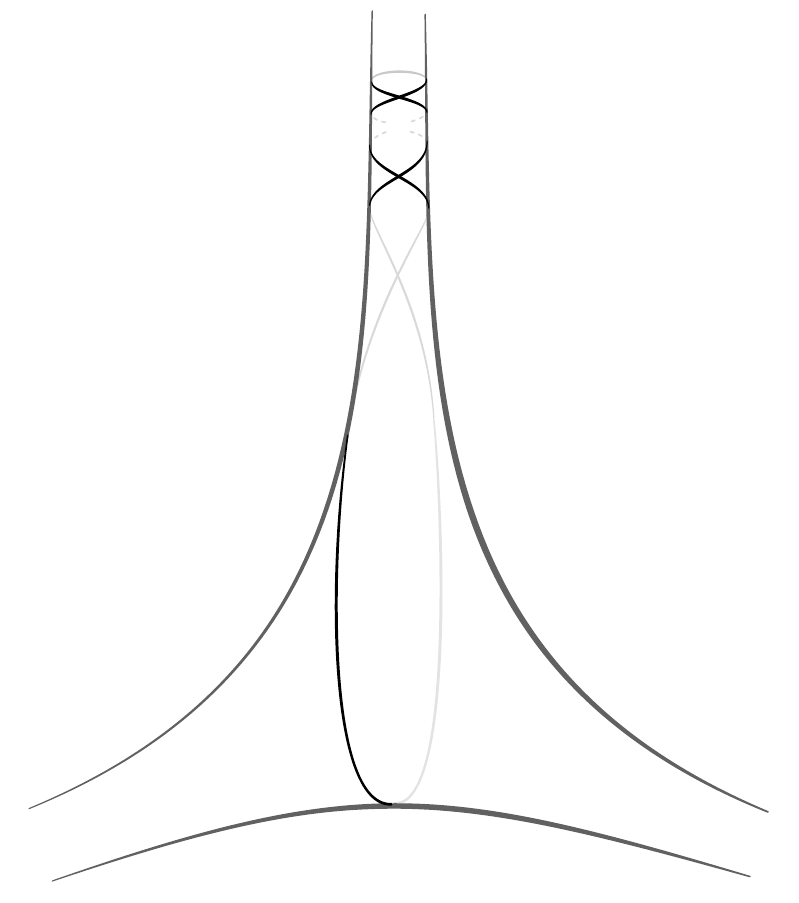}
\vspace{-24pt}
\end{center}
\caption{A corkscrew curve}
\label{fig:Corkscrew}
\end{figure}

By {\it corkscrew} geodesic we mean a geodesic in the homotopy class as described in Figure \ref{fig:Corkscrew}: that is a curve consisting of the concatenation of a simple arc and another that winds $k$-times around a boundary. We can also make the statement "for large enough $k$" explicit: our results hold for all $k \geq 10^{13350}$. As one might imagine, we have not tried to optimize this lower bound.

The problem we solve in the above theorem was first introduced in \cite{Basmajian2013} where such curves were shown to satisfy length bounds, namely upper and lower bounds that are roughly logarithmic in $k$. The upper bound follows from an explicit construction and the corkscrew geodesics above are exactly the ones provided in \cite{Basmajian2013} and they have length exactly $2\arccosh(2k+1)$. They are the natural candidate curves that one might expect as they lie on the smallest possible hyperbolic surface. Our result shows that these candidate curves, at least for sufficiently large $k$, are indeed minimal among all possible geodesics with at least $k$ self-intersections. 

We also point out that, during the writing of this paper, Shen and Wang \cite{ShenWang} posted a paper which shows that the asymptotic growth investigated in \cite{Basmajian2013} is on the order of $2\log(k)$. While our theorem implies theirs, the methods seem quite different.

One might ask what happens for smaller $k$. It is a result of Yamada \cite{Yamada} that the shortest non-simple geodesic taken among all surfaces is precisely on an ideal pair of pants and has $1$ self-intersection (it is a so-called figure $8$ geodesic which is a corkscrew geodesic with a single twist). More generally, Buser \cite{BuserBook} showed that on any hyperbolic surface, the shortest non-simple closed geodesic is always a figure $8$ curve, from which one could also deduce the result for $k=1$. Other than $k=1$, no other explicit values were known before Theorem \ref{thm:main} above.

The proof is somewhat involved and requires multiple, but mostly elementary, ingredients. One of the main new ingredients is to decompose the surface into thick and thin parts when the notion of thick and thin depends on the intersection number. The whole point is to show that, for large enough $k$, the only optimal curve is the corkscrew geodesic on the ideal pair of pants. For this purpose, another essential ingredient is the use of the ideal pair of pants as a companion surface on which we construct a comparison curve.

The paper is structured as follows: after a preliminary section, we provide a sketch of proof which explains the different ingredients. The final section contains all the missing details. 

\noindent {\bf Acknowledgements.}

The authors thank a variety of people they have discussed these ideas with including Binbin Xu and Peter Buser.

\section{Preliminaries}
We consider orientable hyperbolic surfaces, that is surfaces endowed with a complete Riemannian metric of constant curvature $-1$. For our purposes, it suffices to consider finite-type surfaces. In a slight abuse of notation, by boundary we mean the boundary of the Nielsen core of the surface (so the part of the surface where we've removed infinite funnels). 

By curve we mean the continuous image of a circle which is not freely homotopic to a point or a cusp. We assume curves to be primitive, that is not powers of another curve, and not oriented. Generally we think of a curve as a free homotopy class of curve, and we make readily use of the fact that there is a unique geodesic in each free homotopy class. 

The self-intersection of a curve is trickier to define than one might assume. A relatively easy way is to define it as the number of double points (counted with multiplicity) of a geodesic representative. Note that this only works because geodesics minimize intersection and we've assumed that our curves are primitive. We will sometimes refer to a closed geodesic with at least $k$ self-intersection points as a {\it $k$-geodesic}. For instance, all closed geodesics are $0$-geodesics. 

Given a hyperbolic surface $X$ and a point $p\in X$, the injectivity radius of $X$ at $p$ is the radius of the largest open disk embedded at $p$ and will be denoted $\injrad_X(p)$. Given a threshold radius $r>0$, the surface $X$ can be decomposed into its $r$-thin part (resp. thick part), corresponding to the set of points where the injectivity radius is at most (resp. at least) $r$. 

With this in hand, we can now define the quantities we are interested in. Set
$$
M_k:= \inf\{\ell_Y(\alpha) \mid Y \text{ is a hyperbolic surface and }i(\alpha,\alpha)\geq k\}
$$
We will generally denote by $\gamma_k$ or $\gamma$ the curve which we assume is optimal among all geodesics with at least $k$ self-intersections taken among all geodesics on any hyperbolic surface. We denote by $X_k$ or $X$ the optimal surface on which $\gamma$ lives. Observe that here we are assuming an existence result which follows from compacity results of completions of moduli spaces and which is explicitly proven in \cite{Basmajian2013}.

A useful tool, introduced by Thurston in \cite{ThurstonStretch}, is as follows.

\begin{theorem}\label{thm:strip}
A surface $Y$ with a boundary geodesic of positive length can be deformed to a surface with a cusp $X$ and such that for all closed curves $\delta$, the corresponding closed geodesics satisfy $\ell_X(\delta) < \ell_Y(\delta)$.
\end{theorem}
We refer to \cite{Papadopoulos-Theret, Parlier} and \cite{DGK} for other proofs, extensions and applications.

As mentioned before, an essential ingredient in our approach is the following result of Vo \cite{Vo}.
\begin{theorem}\label{thm:vo}
For any $s_0>0$ and $\kappa_0\in \N$ there exists a constant $K=K(s_0,\kappa_0)$ such that for any $k\ge K$, for any surface $X$ with $\sys(X)\geq s_0$ and $\kappa(X)\leq \kappa_0$, the shortest $k$-geodesic is a corkscrew geodesic.
\end{theorem}
Note that Vo uses the term $b(a_c)^k$ geodesics for corkscrew geodesics. By "ideal pair of pants" we mean the unique thrice-punctured sphere (a hyperbolic sphere with three cusps). Note that on an ideal pair of pants, the candidate curves are exactly corkscrew geodesics and have length $2\arccosh(2k+1)$. Furthermore, the above result can be quantified for the ideal pair of pants and a punctured torus with a given lower bound for the systole \cite{Vo}. The result for the punctured torus - and in particular the lower bound on the systole length - is exactly tailored for our needs. 

\begin{corollary}\label{cor:vo} 
On the ideal pair of pants, for any $k\ge10^{125}$, the shortest $k$-geodesic is a corkscrew geodesic with exactly $k$ self-intersections. 
On a punctured torus with systole at least 
$$ 2 \arcsinh\left( \frac{1}{\sinh\left(\frac{5}{2}\right)}\right)
$$
for any $k\ge 10^{13350}$, the shortest $k$-geodesic is a corkscrew geodesic with exactly $k$ self-intersections. 
\end{corollary}

It is these results in particular that we need in Corollary \ref{cor:pants} and Theorem \ref{thm:maintext} at the end of the paper. 

\section{Sketch of proof}
Before proving the needed estimates, we provide a sketch of proof which gives the outline and strategy of our approach. One can imagine that a seasoned hyperbolic geometer would be able to fill in the details themselves, and can then skip the rest of the paper. 

The approach requires dividing the surface into thick and thin parts, where what we mean by thick and thin will be carefully chosen. In what follows, $k$ is always going to be our self-intersection threshold which, at the level of our sketch proof, we can assume is as large as we need it to be at each step. There are finitely many steps. 

Recall that our ultimate goal is to show that for large enough $k$, the only optimal geodesic and corresponding surface are on ideal pairs of pants. These geodesics have length exactly $2 \arccosh(2k+1)$, which is asymptotic to $2\log(k)$ for large $k$.

\noindent{\it The curve $\gamma$ fills $X$}

A very useful technique is to use Theorem \ref{thm:strip} introduced above. It shows that $\gamma$ is a filling curve on $X$, and that $X\setminus \gamma$ is a collection of disks and once-punctured disks where, geometrically, the punctures are realized as cusps. This is because otherwise one can apply a stretch path to strictly reduce the length of all curves (including $\gamma$) which means that $\gamma$ was not optimal. 

\noindent{\it Relating intersection and length in the thick part}

A first key idea, which is already present in most papers relating length to intersection, is to observe that self-intersection of an arc $c$ in the (uniformly) thick part of a surface (where thick is fixed) is quadratically related to its length. That is
$$
\ell(c)^2 \gtrsim i(c,c)
$$
where by $\gtrsim$ we mean bigger than up to a positive multiplicative constant. So on such a surface, its length is roughly bounded below by the square root of its self-intersection number. Observe this immediately prevents, for large enough $k$, a surface which is uniformly thick to be optimal, simply because $\sqrt{k}$ grows a lot faster than $\log(k)$.

This argument can be pushed further to investigate a type of relative thick part, where the thickness depends on $k$. Our measure of thickness will be of the type $\frac{1}{k^a}$ for a (fixed) positive $a$, generally between $0$ and $\frac{1}{2}$. A key step is now to estimate how much length is necessary to go from a portion of a surface which is thick (of say injectivity radius $1$) to the $\frac{1}{k^a}$-thin part.

\noindent{\it Being relatively thin}

Being thin can take two forms: a thin part is either the collar of a short simple closed geodesic or it is the thin part of a horodisk neighborhood of a cusp. In either case, one can estimate that it takes at least $a \log(k) + O(1)$ length to reach the thin part. A precise statement, which works in both situations, is given in Lemma \ref{lem:gotothin}.

Now because we only have "$2\log(k)$" worth of length, for fixed $a$, this means that the curve cannot contain too many arcs that travel from thick to $\frac{1}{k^a}$ thin. In particular, the curve cannot travel to the $\frac{1}{k^a}$ thin part for $a>1$, but we will be using another cut off point. A choice that will work for us will be to use $a=\frac{2}{5}$. It is somewhat arbitrary, but it has the important property and that the curve can contain at most $2$ arcs that travel back and forth between the thick and thin regions, and if there are indeed $2$ of them, the arcs will be forced to travel even less into the very thin regions. By investigating the behavior of arcs in these relatively thin regions, this in turn will allow us to bound the intersection provided by these arcs, but before explaining how this works, we make a crucial observation. 

\noindent{\it Living in the $\frac{1}{k^{\frac{2}{5}}}$-thick part}

For $k$ large enough, we consider the portion of $\gamma$ that lives in the $\frac{1}{k^{\frac{2}{5}}}$-thick part. We denote by $X^T$, where $T$ is for Thick, the $\frac{1}{k^{\frac{2}{5}}}$-thick part of the optimal surface $X$, and by $X^t$, where $t$ stands for thin, the complement of $X^T$. Similarly, we denote $\gamma^T$ the portion of $\gamma$ that lies in $X^T$ and by $\gamma^t$ the portion of $\gamma$ that lies in $X^t$. 

We then show that $\gamma^T$ cannot contribute significantly to intersection. Namely, for large enough $k$, we have 
$$
i(\gamma^T,\gamma^T) < k^{\frac{9}{10}}.
$$
This is a consequence of Lemma \ref{lem:relativethick}. Again, there is a choice involved in the statement. Here we choose $\frac{9}{10}$ because it is bigger than $2$ times $\frac{2}{5}$ but less than $1$. 

For simplicity, we use this technique throughout - that is to choose explicit constants slightly larger than what is strictly needed - to avoid having to deal with logarithmic terms and other lower order terms.

Now this means that $\gamma^t$ must have (at least) self-intersection $k-k^{\frac{9}{10}}$. As we will show, this will force it to enter an even thinner region of the surface.

\noindent{\it Thin collars about geodesics don't help}

We begin by showing that on an optimal surface $X$, 
the $\frac{1}{k^{\frac{2}{5}}}$-thin part $X^t$ does not contain any simple closed geodesics. 

To show this, first observe that if $X^t$ did contain a simple closed geodesic then it would be of length less than $\frac{2}{k^{\frac{2}{5}}}$ and $\gamma$ would have to cross it. In particular $\gamma$ would contain a simple arc that goes from one end of the collar to the other which requires roughly $\frac{4}{5}\log(k)$ of length. An important subtlety here is that we are analyzing the length of the simple arc on $X$, not on $X^t$, although we only analyse intersection in $X^t$. Namely, the natural collar about the putative simple closed geodesic. 

If there are no further arcs that enter the collar in $X^t$, then that portion of $\gamma^t$ does not contribute any intersection and we continue to study the remainder of $\gamma^t$, but with less length at our disposal and we'll reach a contradiction later. So in order for there to be intersection inside the collar, as the initial arc is simple, there must be another arc that enters the collar.

\begin{figure}[h]
\begin{center}
\includegraphics[width=6cm]{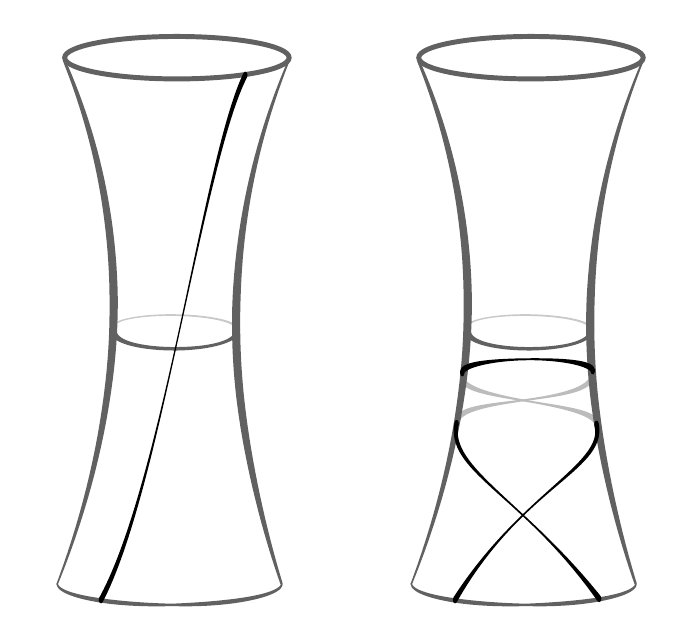}
\vspace{-24pt}
\end{center}
\caption{A transversal arc and a returning arc.}
\label{fig:transversal-returning}
\end{figure}

This second arc is either a simple arc that traverses the collar or a non-simple arc that enters and exits on the same side of the collar. We will refer to these two types of arcs as {\it transversal arcs} or {\it returning arcs} (see Figure \ref{fig:transversal-returning}). 

In both cases, this additional arc has to have length at least roughly $\frac{4}{5}\log(k)$. Note that this prevents a third arc from entering the collar in the thin region. Furthermore, we now have a better lower bound on the length of the core curve in the collar. The shorter the core is, the longer the transversal arc must be. And, by considering the length created by the second arc, the transversal curve can be at most on the order of $\frac{6}{5}\log(k)$. This gives a lower bound on the length of the core curve on the order of $\frac{1}{k^{\frac{3}{5}}}$, a considerable improvement over the previous lower bound which was on the order of $\frac{1}{k}$ (see Lemma \ref{lem:core} for the statement we use). 

We can now look at the two arcs and, using the lower bound on the length of the core curve, give an upper bound on the their intersection and self-intersection (Proposition \ref{prop:coreintersection}). Essentially this boils down to the fact that in a collar, intersection and self-intersection is controlled linearly by the amount of times you wind around the core curve (a type of winding number). As the core curve has length at least on the order of $\frac{1}{k^{\frac{3}{5}}}$ and the total length of $\gamma$ is at most on the order of $2 \log(k)$, this means that intersection is bounded above by a function on the order of $\log(k) k^{\frac{3}{5}}$. Again, this contributes (relatively) little to the needed $k$ self-intersection even if we use up all the length. 

\noindent{\it All the action happens in the cusp}

We can now safely assume that the majority of the intersection happens in the cusp, that is at least $k + o(k)$, and again we analyze the behavior in the cusp. Strands that enter a cusp region are necessarily returning strands and, by virtue of Lemma \ref{lem:gotothin}, if they reach the cusp region of $X^t$, then they have length at least on the order of $\frac{4}{5} \log(k)$. In particular, there can be at most two of them. 

We first consider the case where there are two strands. By performing the same analysis as in the thin collar case, and for essentially the same reasons, we can conclude that the intersection is bounded above by a function on the order of $\log(k) k^{\frac{3}{5}}$. We can thus conclude that there must be a single strand.

The amount of intersection provided by a single strand depends entirely on how high into the cusp the strand travels. As it must provide on the order of $k$ self-intersection points, it must travel to the approximately $\frac{1}{k}$ thin part (see Lemma \ref{lem:gettinghigh}). This requires roughly $2 \log(k)$ length. Note that this estimate shows that (the length of the shortest $k$-geodesic) $M_k$ has $2\log(k)$ asymptotic growth (which is the main result of \cite{ShenWang}). 

\noindent{\it The final blow}

We are now in a position to put serious restrictions on what type of curves the $\gamma_k$ can be. The fact that we have a strand that travels high into the cusp allows us to prove Proposition \ref{prop:onelongarc} which says that outside of the arc that travels high in the cusp, the rest of the curve has uniformly bounded length. In the statement we prove it for $5$, but while we make use of the fact that it is quite short, the actual amount is not that important. 

This proposition has an essential corollary: we can deduce that $\gamma_k$ lies in an ideal pair of pants (Corollary \ref{cor:pants}). The corollary is proved using the isoperimetric inequality (in the hyperbolic plane and the once-punctured hyperbolic plane).

With this in hand, we are now in position to apply Vo's result (Corollary \ref{cor:vo}). It states that for all $k>10^{125}$, the shortest $k$-geodesics on an ideal pair of pants are exactly our candidate curves and this completes the proof. 

We now provide the details necessary for the above steps. 

 \section{The structure of optimal geodesics for large enough $k$}
 
One of the key aspects of our approach is to use the fact that it costs length to reach the thin part. This is quantified by the following lemma. 
 
 \begin{lemma}\label{lem:gotothin}
Fix a constant $N>e$ and a hyperbolic surface $X$.

An arc $c$ that goes from $p$ to $q$ with $\injrad_X(p)= \arcsinh( 1)$ and $\injrad_X(q) \leq \frac{1}{N}$ satisfies
$$
\ell(c) > \log(N) - 1.
$$
 \end{lemma}
 
 \begin{proof}
 The point $q$ lies in the half-collar of a curve of length at most $\frac{2}{N}$ (this length could be $0$ in which case it lies in a horocyle neighborhood of a cusp). 
 
 If $p$ doesn't lie in the same half-collar then any arc joining $p$ and $q$ must pass through another point $p'$ with injectivity radius $\arcsinh(1)$ so we may as well assume that $p$ lies in the same half-collar as $q$. In this case the loops $\gamma_p$ and $\gamma_q$ of lengths $2 \arcsinh(1)$ and less than 
 $\frac{2}{N}$, based respectively at $p$ and $q$, are isotopic. In this isotopy class there is a core closed geodesic $\alpha$ (or a cusp which we can think of as a core geodesic of length $0$).
We may assume that the length of $\gamma_q$ is exactly 
$\frac{2}{N}$ since any geodesic loop whose length is less than that is even further from $p$.

We consider the cylinder $C$ bounded by the loop $\gamma_p$ based in $p$ on one end and by $\alpha$ on the other. Consider $\delta$ the unique geodesic of shortest length between $p$ and $\alpha$ (which is of infinite length if $\alpha$ is a cusp). 

The point $q$ belongs to the circle of points of this cylinder that are the base points of geodesics loops of length $\frac{2}{N}$. This circle of points is parallel to $\alpha$ (or is a horocycle in the case of a cusp). Of all points of this circle of points, the closest one to $p$ is exactly the intersection between this circle and $\delta$. Hence, as we want a lower bound on the distance between $p$ and $q$, we can assume that $q$ is exactly this point (see Figure \ref{fig:halfcylinder}). 

\begin{figure}[h]
\leavevmode \SetLabels
\L(.51*.78) $p$\\%
\L(.5*.42) $q$\\%
\L(.505*.54) $\delta$\\%
\L(.49*.00) $\alpha$\\%
\endSetLabels
\begin{center}
\AffixLabels{\centerline{\includegraphics[width=3.5cm]{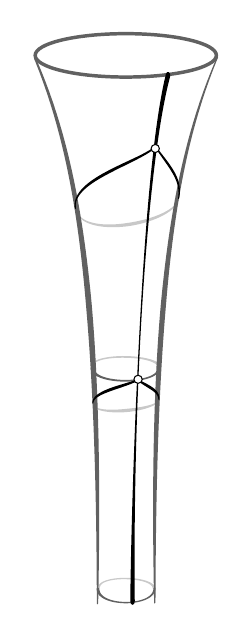}}}
\vspace{-24pt}
\end{center}
\caption{The distance between $p$ and $q$ is minimal when they both lie on $\delta$}
\label{fig:halfcylinder}
\end{figure}
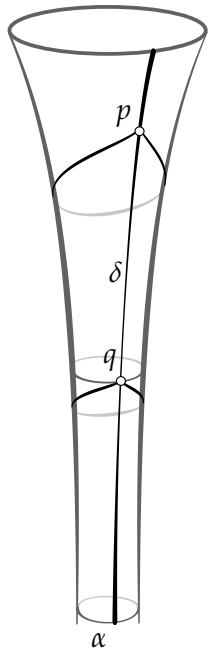

If we cut the cylinder $C$ along $\delta$, we obtain a quadrilateral or, in the case where we have a cusp, a triangle with an ideal vertex. It has an axis of reflection, say $\beta$, which cuts the quadrilateral into two quadrilaterals with two right angles as in Figure \ref{fig:Quad-reflect}. 

\begin{figure}[h]
\leavevmode \SetLabels
\L(.195*.9) $p$\\%
\L(.193*.08) $p$\\%
\L(.51*.79) $q$\\%
\L(.51*.15) $q$\\%
\L(.35*.83) $t$\\%
\L(.107*.61) $\arcsinh(1)$\\%
\L(.48*.56) $\sfrac{1}{N}$\\%
\endSetLabels
\begin{center}
\AffixLabels{\centerline{\includegraphics[width=10cm]{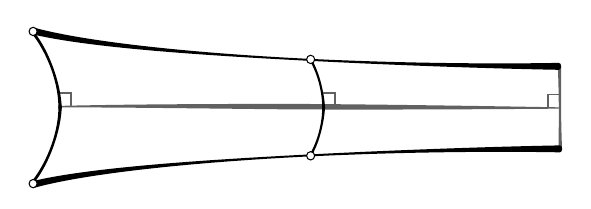}}}
\vspace{-24pt}
\end{center}
\caption{Quadrilaterals}
\label{fig:Quad-reflect}
\end{figure}
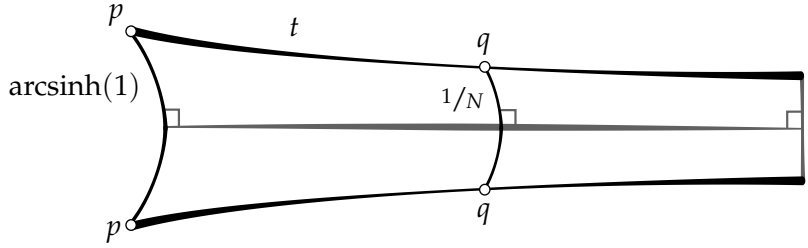

As $q$ lies on $\delta$, we can represent the loop $\gamma_q$ on the quadrilateral and it too is invariant by the reflection. We denote by $t$ the distance between $p$ and $q$ on the quadrilateral (and on the surface). 

Our goal is to compute a lower bound for $t$ using elementary hyperbolic geometry. Our first argument is geometric: as the lengths $\arcsinh(1)$ and $\frac{1}{N}$ are fixed, the geometry of the quadrilaterals, and hence $t$, is entirely determined by the length of the core geodesic $\alpha$. Observe that $t$ behaves monotonically in the length of $\alpha$, and is shortest when $\ell(\alpha)=0$. This can be seen geometrically - the arc of length $\frac{2}{N}$ has to move closer to $p$ as the length of $\alpha$ decreases - or by a standard computation. In any event, as we are interested in a lower bound on $t$, we can assume the length of $\alpha$ is $0$ and we are in the case where the quadrilateral has degenerated to an ideal triangle as in Figure \ref{fig:Idealtriangle}. 

It is fortuitous that in this case, the computations become easier. Using standard hyperbolic trigonometry, the (half)-angle at $p$ is $\frac{\pi}{4}$ as in Figure \ref{fig:Idealtriangle}.

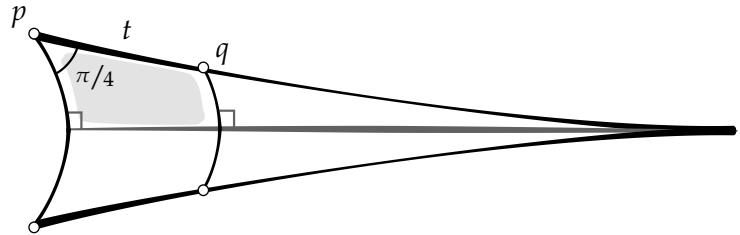
\begin{figure}[h]
\leavevmode \SetLabels
\L(.195*.9) $p$\\%
\L(.37*.76) $q$\\%
\L(.25*.65) $\sfrac{\pi}{4}$\\%
\L(.29*.83) $t$\\%
\endSetLabels
\begin{center}
\AffixLabels{\centerline{\includegraphics[width=10cm]{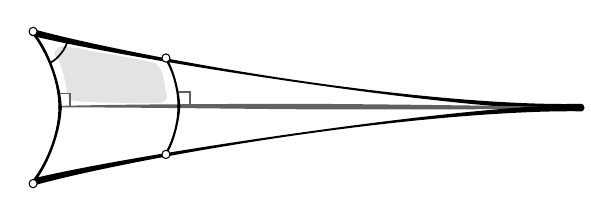}}}
\vspace{-24pt}
\end{center}
\caption{The worst case scenario is a triangle with an ideal vertex}
\label{fig:Idealtriangle}
\end{figure}
Now we can compute $t$ using the quadrilateral with two right angles and sides of length $\arcsinh(1)$, $\frac{1}{N}$, $t$ and angle at $p$ equal to $\frac{\pi}{4}$ (shaded in Figure \ref{fig:Idealtriangle}). Using the formula for such quadrilaterals (see for instance \cite[Section VI.3, pp. 88-89]{Fenchel}) we have
$$
\sinh\left(\frac{1}{N}\right)= \sinh(\arcsinh(1)) \cosh(t) - \cosh(\arcsinh(1)) \sinh(t) \cos\left(\frac{\pi}{4}\right)
$$
and thus
$$
\sinh\left(\frac{1}{N}\right)= \cosh(t) - \sinh(t) = e^t.
$$
From this we obtain 
$$
t= \log\left(\frac{1}{\sinh\left(\frac{1}{N}\right)}\right)>\log(N) - 1
$$
for all $N>e$ (the last inequality can be verified by a quick calculus computation). This completes the proof. 
\end{proof}
 
 Note that the constant $1$ in the above statement is chosen for simplicity. The quantity
 $$ \log\left(\frac{1}{\sinh\left(\frac{1}{N}\right)}\right) -\log(N)$$
 goes to $0$ as $N$ increases. The lower bound on $N$ ensures that the loop around $q$ is inside the half-collar bounded by the loop around $p$. Of course it is chosen so that the quantity $\log(N)-1$ is positive. 
 
For a real number $\frac{1}{2}>a>0$ we provide a general upper bound on the amount of intersection that can occur in the $\frac{1}{k^a}$-thick part of $X$. More precisely, consider $\gamma_a:= \gamma \cap X^{\frac{1}{k^a}{\text -} T}$. 

 \begin{lemma}\label{lem:relativethick}
$ i(\gamma_a,\gamma_a) \leq (2 k^a \log(k))^2 $
\end{lemma}
 
 Note that the statement is true for all $a>0$ but as by definition the self-intersection of $\gamma$ is $k$, for $a\geq \frac{1}{2}$ the lemma is void of content. 
 
 \begin{proof}
 We break up $\gamma_a$ into a maximal collection of segments of length at most $\frac{1}{k^a}$. Observe that any two of these segments can intersect at most once. 
 
 As $\ell(\gamma_a)\leq \ell(\gamma)\leq 2\arccosh(1+2k) <2 \log(k) + 4$, the maximal number of segments is
 \begin{eqnarray*}
 N \leq & \left[\frac{\ell(\gamma_a)}{\sfrac{1}{k^a}} \right] 
 <k^a (2 \log(k) +4) + 1 <2 \sqrt{2} k^a \log(k)
\end{eqnarray*}
where the last inequality holds for all $k> 10^6$. 

And so the total number of intersections is $N(N-1)/2<N^2/2$ and hence at most
$$
\frac{1}{2}(2\sqrt{2} k^a \log(k))^2 = (2 k^a \log(k))^2
$$
as required.
 \end{proof}
 
 Recall that the optimal $k$-geodesic $\gamma$ is filling on the optimal surface $X$, that $\gamma^{T}$ denotes the portion of $\gamma$ in the $\frac{1}{k^{2/5}}$ thick part of $X$, and that we may assume there exists a point in $X$ for which the injectivity radius is at most 
 $\frac{1}{k^{2/5}}$.

\begin{corollary}\label{cor:intersectionboundinthick}
For large enough $k$, $i(\gamma^T,\gamma^T) < k^{\frac{9}{10}}.$
\end{corollary}

\begin{proof}
Setting $a=\frac{2}{5}$ in lemma \ref{lem:relativethick},
it suffices to observe that for large enough $k$, 
$$(2 k^{\frac{2}{5}} \log(k))^2 <k^{\frac{9}{10}}.$$
It boils down to checking for which $k$ we have 
$$
4 (\log(k))^2 < k^{\frac{1}{10}}
$$
which holds for all $k\geq 10^{48}$. 
\end{proof}

We now turn to the part of $\gamma$ in the $\frac{1}{k^{\frac{2}{5}}}$-thin part denoted $X^t$.
The following lemma shows that around collars of curves in $X^t$, the behavior of $\gamma$ is well-controlled.

\begin{lemma}\label{lem:thinarcs}
Suppose $k$ is large. Let $\sigma$ be a simple closed geodesic in $X^t$ and $C_\sigma$ be the collar region of $\sigma$ in $X^t$. Then $\gamma \cap C_\sigma$ consists of at most two arcs.
\end{lemma}

\begin{proof} This is more or less a direct consequence of Lemma \ref{lem:gotothin}, where $N=k^{\frac{2}{5}}$.

Indeed all strands inside $C_\sigma$ must reach and leave the $\frac{1}{k^{\sfrac{2}{5}}}$-thin region. There is a somewhat subtle point here: we're looking at the length of the extension of the arc that enters the thin part of the collar. These full arcs, by Lemma \ref{lem:gotothin} above, are of length at least 
$$
 2\left( \frac{2}{5} \log(k) -1 \right).
$$
Now if there were at least $3$ then the total length of $\gamma$ would be at least
$$
3 \left( \frac{4}{5} \log(k) -2 \right) = \frac{12}{5} \log(k) -6 > 2 \log(k) + 4 > \ell(\gamma)
$$
which provides a contradiction if $\log(k) > 25$, that is if $k>e^{25}$. 
\end{proof}

Using the lemma above, we can now show that if there is intersection in the collar of a curve in $X^t$, then the core curve cannot be too short. 

\begin{lemma}\label{lem:core}
Suppose $k$ is large. Let $\sigma$ be a simple closed geodesic in $X^t$. Denote by $C_\sigma$ the collar region of $\sigma$ in $X^t$, and by $\gamma_\sigma$ the restriction of $\gamma$ to $C_\sigma$. If $i(\gamma_{\sigma},\gamma_{\sigma})>0$, then 
$$\ell(\sigma) > \frac{1}{k^{\sfrac{4}{5}}}.
$$
\end{lemma}
\begin{proof}
As $\gamma$ is filling, there must be at least one transversal strand in $C_\sigma$. Hence any intersection must come from the intersection between the transversal strand and the returning strand, or from the intersection between two transversal strands. As we suppose that $i(\gamma_{\sigma},\gamma_{\sigma})>0$, there are at least 2 strands (and by the previous lemma, there are exactly $2$).

As each strand has length at least $\frac{4}{5} \log(k) -2$, this means that a transversal strand has length at most 
$$
2 \log(k) + 4 - \left( \frac{4}{5} \log(k) -2 \right) = \frac{6}{5} \log(k) +6.
$$
Now if the core geodesic is of length less than $\frac{1}{k^a}$ then by Lemma \ref{lem:gotothin}, the length of the transversal strand is at least $2 \log(k^a) -2$. Hence for any $a>\frac{3}{5}$, we would reach a contradiction for large enough $k$. We choose $a = \frac{4}{5}$ for simplicity and argue by contradiction. We have 
$$
\frac{6}{5} \log(k) +6 > \frac{8}{5} \log(k) - 2
$$
that is
$$
8 > \frac{2}{5} \log(k)
$$
which fails for $k > e^{20}$ and hence $\ell(\sigma)>\frac{1}{k^\frac{4}{5}}$.
\end{proof}

We now bound the total intersection that occurs in collars of simple closed geodesics inside $X^t$. We denote by $\sigma_1,\hdots,\sigma_n$ the collection of simple closed geodesics of $X^t$ and by $C$ the collection of their collars (so $C = \cup_{i=1}^n C_{\sigma_i}$). We denote by $\gamma_C$ the restriction of $\gamma$ to $C$ (so $\gamma_C= \cup_{i=1}^n C_{\gamma_i}$).

\begin{proposition}\label{prop:coreintersection}
For large enough $k$, there are at most 2 simple closed geodesics in $X^{t}$.
Moreover:
\begin{itemize}
 \item if there are none, then $X$ must contain a cusp.
 \item if there is one,
 then $i(\gamma_C,\gamma_C) < k^{\frac{9}{10}},$
 \item if there are two, then $i(\gamma_C,\gamma_C)=0.$
\end{itemize}

\end{proposition}


\begin{proof}
Note that there are at most two strands in total that lie in $C$, and as $\gamma$ must fill there are at most two simple closed geodesics in $X^t.$ 

Suppose $X^t$ has no simple closed geodesics. If $X$ does not contain a cusp neightborhood, then $X=X^T$, $\gamma=\gamma^{T}$, and hence by Corollary \ref{cor:intersectionboundinthick}, $i(\gamma,\gamma) < k^{\frac{9}{10}},$ leading to a contradiction since 
$i(\gamma,\gamma)\geq k$.

If $X^t$ contains two simple closed geodesics then there must be two transversal strands to each of these geodesics, and hence $i(\gamma_C,\gamma_C)=0.$

We are left to prove the case where there is one simple closed geodesic in $X^t$. Again if there is one transversal then 
$i(\gamma_C,\gamma_C)=0.$ Thus the only non-trivial case to consider is when there is one simple closed geodesic in $X^t$ and two strands, one of which must be transversal. There are cases to consider:

\noindent{\it Case 1:} Both strands are transversal.

Denote the two strands $a$ and $b$ and assume that $i(a,b)>1$. Note that both are individually simple and go from one end of the collar to the other "monotonically": that is, if we parametrized the collar with respect to the core geodesic and distance to the core geodesic, this distance would be monotonic. The main observation is that two consecutive intersection points are joined by one arc of $a$ and one arc of $b$. These are not homotopic arcs, and their concatenation is freely homotopic to the core geodesic. Hence it follows that
$$
(i(a,b)-1) \ell(\sigma) < \ell(a) + \ell(b) 
$$
(Note that we need at least $2$ intersection points for this to work and each extra intersection point adds a new copy of a curve homotopic to $\sigma$.)

As $\ell(a) + \ell(b) < \ell(\gamma)< 2 \log(k) + 4$ and by lemma 
\ref{lem:core},
$\ell(\sigma) > \frac{1}{k^{\sfrac{4}{5}}}$, we have 

$$
i(a,b) < 1 + k^{\frac{4}{5}} (2 \log(k) + 4)< k^{\frac{9}{10}}.
$$
The last inequality holds for large enough $k$. \\

\noindent{\it Case 2:} There is one transversal and one returning strand.

Let $a$ be the transversal strand and $b$ the returning strand.
The proof works almost identically as in case (1) but with two small changes. 

The first one is that we have to account for the self-intersection of $b$. The second one is that $b$ does not pass through the collar monotonically, but first approaches the core curve before turning back. To overcome the second difficulty, we split $b$ into two simple and monotonic arcs $b^+$ and $b^-$. 

Now we have $i(b,b) = i(b^+,b^-)$. As before, between consecutive intersection points, we have arcs of $b^+$ and $b^-$ that can be concatenated to be homotopic to $\sigma$. Hence 
$$
(i(b,b)-1) \ell(\sigma) < \ell(b^+) + \ell(b^-) = \ell(b).
$$
With the same argument, we can bound the intersection between the arc $a$ and $b^+$, resp. $b^-$: 
$$
(i(a,b^+)-1) \ell(\sigma) < \ell(a) + \ell(b^+). 
$$
and
$$
(i(a,b^-)-1) \ell(\sigma) < \ell(a) + \ell(b^-). 
$$
so 
$$
i(a,b) \leq 2 + \frac{\ell(a) + \ell(b)}{\ell(\sigma)}.
$$
As before, we use the upper bound on the length of $\gamma$ and the lower bound on the length of $\sigma$ to obtain:
$$
i(b,b) < 1 + k^{\frac{4}{5}} (2 \log(k) + 4)
$$
and
$$
i(a,b) < 2 + k^{\frac{4}{5}} (2 \log(k) + 4)
$$
Hence 

$$
i(\gamma_C,\gamma_C) = i(a,b) + i(b,b) < 3 +2 k^{\frac{4}{5}} (2 \log(k) + 4)<k^{\frac{9}{10}}
$$
where the last inequality hold for all $k>10^4$. 
\end{proof}

With these results in hand, for large enough $k$, observe that outside of cusp regions, $\gamma$ self-intersects at most $2 k^{\frac{9}{10}}$ times. Hence there must be $k-2 k^{\frac{9}{10}}$ intersections in cusp regions. In particular, 

\begin{corollary}
For large enough $k$, the optimal surface $X$ must have a cusp. 
\end{corollary}

\section{Optimal excursions into cusped surfaces and the main result}

We now prove a general estimate that relates length to intersection in a cusp region. Similar estimates can be found in the literature \cite{Basmajian2015} but here we follow \cite{Vo}. 


\begin{lemma}\label{lem:intercusp}
Let $a$ be a returning strand that enters and exits a cusp region at points on the horocycle boundary of length $1$ and whose highest point has injectivity radius $\frac{1}{N}$ for some $N>2$. Then 
\[
2\sqrt{\frac{1}{4\sinh^2\frac{1}{N}}-1}-2
\le i(a,a)
\le 2\sqrt{\frac{1}{4\sinh^2\frac{1}{N}}-1}-1
\le N-1.
\]
Similarly, if $a$ is as above and $b$ is another one which is not shorter than $a$, then 
\[
i(a,b) 
\le 4\sqrt{\frac{1}{4\sinh^2\frac{1}{N}}-1}
\le 2N.
\] 
\end{lemma}
\begin{proof}
Upon lifting to the upper-half plane model for the hyperbolic plane, we can normalize so that the horocycle boundary of length $1$ lifts to the Euclidean line $y=1$ and the parabolic element identified with the element in the fundamental group that wraps once around the cusp is $f:z\mapsto z+1$. The lift of $a$ is a geodesic segment $\tilde{a}$ that is contained on a Euclidean semicircle orthogonal to the real axis. This semicircle can be further normalized so that its center is the origin. Denote the endpoints of this semicircle by $-r$ and $r$. A computation yields 
\[
r=\cosh \frac{\ell(a)}{2}.
\]
The $f$-translates of $\tilde{a}$ intersects $\tilde{a}$ exactly $i(a,a)$ times. We have
\[
i(a,a)+1=\max\left\{n\in\mathbb{N} : -r+n\le -r+2\sqrt{r^2-1}\right\}
=\floor*{2\sqrt{r^2-1}}. 
\]
Note that a point on the cusp region that has injectivity radius $x$ belong to the horocycle of length $2\sinh(x)$. 
Since the highest point of $a$ has injectivity radius $\frac{1}{N}$, we have
\[
r=\frac{1}{2\sinh\frac{1}{N}}.
\]
Thus
\[
2\sqrt{\frac{1}{4\sinh^2\frac{1}{N}}-1}-2
\le i(a,a)
\le 2\sqrt{\frac{1}{4\sinh^2\frac{1}{N}}-1}-1
\le N-1,
\]
where the last inequality follows from $\sinh(x)>x$ for $x>0$.

For two returning strands $a$ and $b$, where $a$ is as above, and $b$ is not shorter, we have
\[ 
i(a,b) \le 2 \ceil*{w(a)} \le 2(i(a,a)+1),
\]
where $w(\cdot)$ denotes the winding number and the first inequality follows from \cite[Lemma 3.2]{Erlandsson-Parlier}. And so using the above inequality we have
\[
i(a,b) \leq 4\sqrt{\frac{1}{4\sinh^2\frac{1}{N}}-1}\le 2N.
\]
\end{proof}

We now show that $\gamma$ will have at most $1$ strand in cusp regions. 
\begin{lemma}
There is at most one strand of $\gamma$ that enters a cusp neighborhood of $X^t$. 
\end{lemma}
\begin{proof}
The arguments are similar to the earlier arguments regarding strands of $\gamma$ entering collar neighborhoods of the thin part. As before, using Lemma \ref{lem:gotothin} there are at most 2 strands in total in $X^t$, so in particular in a cusp region of $X^t$. Now if there are exactly 2, then neither can go too high. More precisely, each strand has length at least $\frac{4}{5} \log(k) -2$, and this means that the other strand has length at most 
$$
2 \log(k) + 4 - \left( \frac{4}{5} \log(k) -2 \right) < \frac{6}{5} \log(k) +6.
$$
Hence, again by Lemma \ref{lem:gotothin}, if it enters the $\frac{1}{k^a}$-thin region, then it is of length at least 
$$
2 a \log(k) - 2.
$$
So if $a\geq \frac{4}{5}$ then 
$$
 \frac{8}{5}\log(k) -2 < \frac{6}{5} \log(k) +6
$$
and thus
$$
 \frac{1}{5}\log(k) < 4
$$
which is false for all $k>e^{20}$. Thus, for large enough $k$, the arcs never reach a point of injectivity radius $\frac{1}{k^{\sfrac{4}{5}}}$. 
We can now apply Lemma \ref{lem:intercusp}, 
their intersection number is bounded from above by $2k^{\sfrac{4}{5}}$. 
Hence, if there are 2 strands, for large enough $k$ they cannot provide the needed intersection.
\end{proof}

Now that we know that we have a single strand in the cusp region, and that it must contribute at least $k- 2 k^{\frac{9}{10}}$ to intersection, we have the following.

\begin{lemma}\label{lem:gettinghigh}
For large enough $k$, the curve $\gamma$ contains a strand that reaches the 
$\frac{1}{k - 2 k^{\frac{9}{10}}}$
thin part of a cusp region. 
\end{lemma}
 \begin{proof}
 As before, we denote the strand by $a$ and the minimal injectivity radius of the strand $\frac{1}{N}$. We can again apply Lemma \ref{lem:intercusp}, we have
\[ 
i(a,a) < N
\]
 We know that $i(a,a) \geq k - 2 k^{\frac{9}{10}}$ which proves the lemma.
 \end{proof}
 
 At this point we have shown that the optimal $k$-geodesic $\gamma$ (for $k$ large) enters a single cusp in the $X^{t}$ part of the optimal $k$-surface $X$ and must reach the $\frac{1}{k - 2 k^{\frac{9}{10}}}$ thin part of this cusp. In addition, there are more than 
 ${k - 2 k^{\frac{9}{10}}}$ intersections in this cusp part.

We now show that outside of the arc $a$ which we've described above, the rest of the curve is of length at most $5$. The choice of $5$ is somewhat arbitrary and we're not trying to optimize it, but it is small enough to suit our purposes.

\begin{proposition}\label{prop:onelongarc} 
For large enough $k$, the curve $\gamma$ consists of a single arc $a$ which travels into the cusp and a second arc $b$ of length $\ell(b)<5$
\end{proposition}

\begin{proof}
We argue by contradiction. Suppose that the complementary arc of $a$ is of length at least $5$. We shall construct a geodesic $\gamma'$ on the thrice punctured sphere of length shorter than $\gamma$ and with at least as many intersection points. 

Consider the arc $a^t$, the restriction of $a$ to $X^t$. For large enough $k$, its complementary arc (that is $\gamma \setminus a^t$ has at most $2 k^{\frac{9}{10}}$ self-intersections: at most $k^{\frac{9}{10}}$ in $X^T$ (Corollary \ref{cor:intersectionboundinthick}) and at most $k^{\frac{9}{10}}$ in the collars around curves in the thin part (Proposition \ref{prop:coreintersection}). So $a^t$ (and thus $a$) must have self-intersection at least $k - 2 k^{\frac{9}{10}}$. 

\begin{figure}[h]
\leavevmode \SetLabels
\L(.47*.45) $a'$\\%
\L(.548*.2) $h$\\%
\endSetLabels
\begin{center}
\AffixLabels{\centerline{\includegraphics[width=8cm]{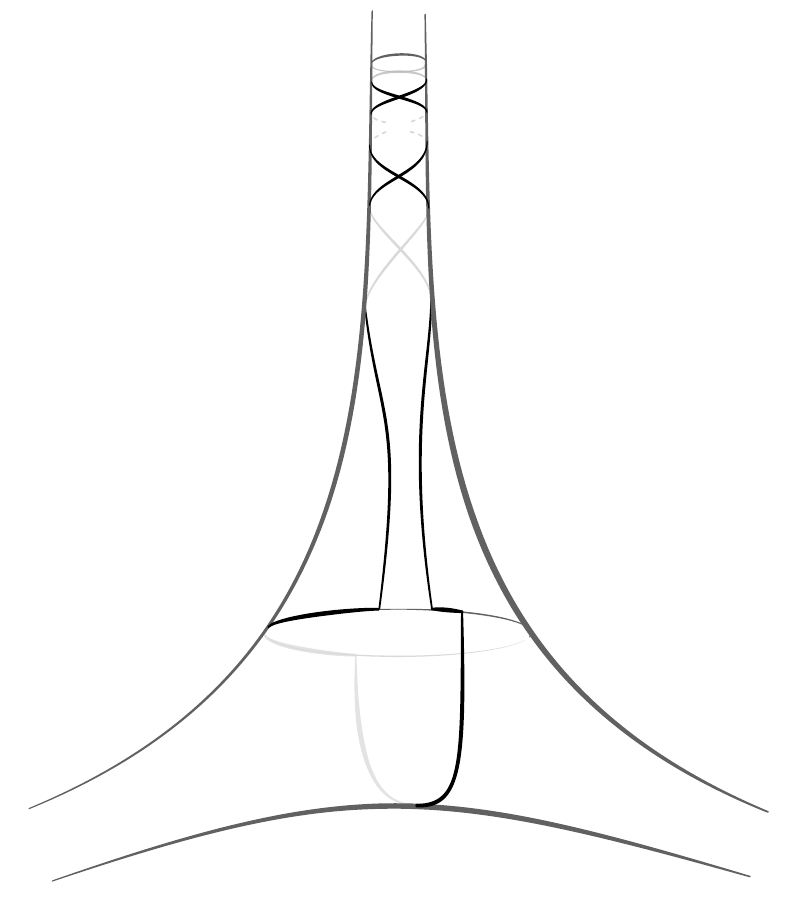}}}
\vspace{-24pt}
\end{center}
\caption{Constructing $\gamma'$}
\label{fig:Gammaprime}
\end{figure}

We now construct the geodesic $\gamma'$ on the thrice punctured sphere as follows. We begin by taking a copy of $a$ on a cusp of our choosing. Note that there are many ways to take a copy of $a$ (by rotating any given copy) but we don't mind which one we choose here. We add $2 k^{\frac{9}{10}}+1$ intersection points to $a$ to obtain an arc $a'$ by adding extra winding by following the shortest geodesic loop based at the highest point. 
By Lemma \ref{lem:gettinghigh}, this shortest geodesic loop has length 
\[
L \le \frac{2}{k - 2 k^{\frac{9}{10}}}
< \frac{1}{2 k^{\frac{9}{10}}+1}
\]
for $k>10^{10}$, 
so the length of $a'$ is at most 
\[
L (2 k^{\frac{9}{10}}+1) + \ell(a)
< \ell(a) + 1.
\]
We now join the endpoints of $a'$ in the following way. We use the shortest orthogeodesic $h$ on the ideal pair of pants with endpoints on the horocycle of length $1$, and add paths that follow the horocyle boundary and remain disjoint (see Figure \ref{fig:Gammaprime}). This extra arc is of length at most $1$ plus the length of $h$ which by a routine calculation on the ideal pants is at most $4\log(2)$.
So the total length of $\gamma'$ is strictly less than
$\ell(a)+1+1+4\log(2)<\ell(a)+5.$ 
Moreover, $\gamma'$ has (at least) one extra intersection point compared to $\gamma$ which we have supposed is of length $\geq 5 + \ell(a)$.
 
Now, though the curve we are considering is not geodesic, its homotopy class is clear, and so we control the self-intersection of the curve very well here. Finally, its geodesic representative will have length less than that of $\gamma$, leading to a contradiction.
\end{proof}

\begin{corollary}\label{cor:pants}
For large enough $k$, the surface $X$ is a pair of pants.
\end{corollary}
\begin{proof}
We use the isoperimetric inequalities of the hyperbolic plane and the once-punctured hyperbolic plane, see, e.g., \cite{Adams-Morgan1999}. The former says that a rectifiable simple curve of length $L$ that surrounds a region of area $A$ in the hyperbolic plane satisfies
$$
A^2 + 4\pi A \leq L^2.
$$
The latter says that a rectifiable simple curve of length $L$ that surrounds a once-punctured (that is once-cusped) hyperbolic region of area $A$ satisfies
$$
A \leq L.
$$
Note that the latter inequality is implied by the former. 

Now we cut along the arc $b$ which is of length $<5$ by Proposition \ref{prop:onelongarc}, and along the horocycle of length $1$ around the cusp. The result is a collection of disks or once-punctured disks. As each curve contributes boundary length to exactly two regions, the total length $\ell$ of the boundaries of the disks is strictly less than $2(5+1)=12$. Now for each disk, punctured or not, if $\ell_i$ is the boundary length and its area is $a_i$, then it satisfies $a_i \leq \ell_i$. By adding them up we have that the area of the surface satisfies
$$
{\rm Area} (X) \leq \ell < 12.
$$
Now 
$$
{\rm Area} (X) = - 2\pi \chi(X) <12
$$
where $\chi(X)$ is the Euler characteristic of $X$, and hence 
$$
-\chi(X) = 2g-2+n < \frac{6}{\pi}<2
$$
and as it is an integer we have that $-\chi(X)=1$ and $X$ is either a punctured torus or a pair of pants. 
Suppose $X$ is a punctured torus. As $b$ fills, it must cross the systole of $X$ and its collar neighborhood. Now as $b$ is of length less than $5$, a computation using the collar lemma, gives us a lower bound on the systole of $X$ as follows. If $\ell_X$ is the systole length then 
$$
\ell(b)>2\arcsinh\left( \frac{1}{\sinh\left(\frac{\ell_X}{2}\right)}\right).
$$
Hence
$$
\ell_X > 2 \arcsinh\left( \frac{1}{\sinh\left(\frac{5}{2}\right)}\right).
$$
Now by Corollary \ref{cor:vo}, for $k\ge 10^{13350}$, the shortest curve with at least $k$ self-intersections is a corkscrew geodesic and thus cannot fill $X$, leading to a contradiction. Thus $X$ must be a pair of pants and moreover, since it only has only cusps as ends, it is an ideal pair of pants. 
\end{proof}

We can now deduce the main theorem. We recall its statement.

\begin{theorem}\label{thm:maintext}
Let $\gamma_k$ be the closed geodesic on a connected hyperbolic surface $X_k$ of minimal length among all primitive closed geodesics with at least $k$ self-intersection points taken among all closed geodesics on all possible hyperbolic surfaces. Then for large enough $k$, $X_k$ is an ideal pair of pants and $\gamma_k$ is a corkscrew geodesic with exactly $k$ self-intersections.
\end{theorem}

\begin{proof}
By the above (Corollary \ref{cor:pants}), $X=X_k$ is an ideal pair of pants. Now by Vo's theorem \cite{Vo}, for all $k\ge 10^{13350} > 10^{125}$, the shortest geodesic on an ideal pair of pants with at least $k$ self-intersections is necessarily a corkscrew geodesic with exactly $k$ self-intersections and we are done. 
\end{proof}

{\it Addresses:}\\
The Graduate Center, City University of New York, NY, New York, 10016 and Hunter College, City University of New York, 695 Park Ave., New York, NY, 10065\\
DMATH, FSTM, University of Luxembourg, Esch-sur-Alzette, Luxembourg\\
School of Mathematical and Statistical Sciences, Arizona State University, USA

{\it Emails:}\\
abasmajian@gc.cuny.edu\\
hugo.parlier@uni.lu\\
thihanhv@asu.edu

\end{document}